\documentclass[a4paper,12pt]{amsart}
\usepackage[letterpaper, margin=1.2in]{geometry}
\usepackage{latexsym}
\usepackage{amssymb}
\usepackage{amsmath}
\usepackage{amsbsy}
\usepackage{amsfonts}
\usepackage[table]{xcolor}
\usepackage{pgfplots}
\usepackage{color}
\usepackage{multirow}
\usepackage{graphicx}
\usepackage{txfonts,pxfonts} 
\usepackage{tikz}
\usetikzlibrary{calc,arrows}
\usetikzlibrary{calc}
\usepackage{verbatim}
\everymath{\displaystyle}
\newtheorem{thm}{Theorem}[section]

\newtheorem{lem}[thm]{Lemma}

\theoremstyle{definition}

\newtheorem{rem}{Remark}

\newcommand{\Z}{\mathbb Z}

\def\ol{\overline}

\def\la{\lambda}

\def\md#1{\ \mbox{\rm(mod }{#1})}

\newcounter{cs}
\stepcounter{cs}
\newcommand{\casos}{\begin{itemize}}
	\newcommand{\fcasos}{\end{itemize}\setcounter{cs}{1}}

\newfont{\tit}{cmr12 scaled \magstep3}
\title[The index divisors and monogenity of certain nonic number fields]{On the index divisors and monogenity of certain nonic number fields}
\author[Omar Kchit]{Omar Kchit}
\address{Faculty of Sciences Dhar El Mahraz, P.O. Box  1796 Atlas-Fes, Sidi Mohamed ben Abdellah University,  Morocco}
\email{omar.kchit@usmba.ac.ma, \,\, orcid: 0000.0002.0844.5034}
\begin{document}
	\keywords{Theorem of Dedekind,  Theorem of Ore, prime ideal factorization,  Newton polygon, Index of a number field, Power integral basis, Monogenic}
	\subjclass[2010]{11R04, 11Y40, 11R21}
	
	\begin{abstract}
		In this	paper, for any nonic number field $K$ generated by a root $\alpha$ of a monic irreducible trinomial $F(x)=x^9+ax+b \in \mathbb{Z}[x]$ and for every rational prime $p$, we characterize when $p$ divides the index of $K$. We also describe the prime power decomposition of the index $i(K)$. In such a way we give a partial answer of Problem $22$ of Narkiewicz (\cite{Nar}) for this family of number fields. In particular if $i(K)\neq 1$, then $K$ is not mongenic. We illustrate our results by some computational examples. 
	\end{abstract}
	\maketitle
	\section{Introduction}
	Let $K$ be a number field of degree $n$ and $\mathbb{Z}_K$ its ring of integers. For any primitive element $\alpha\in \mathbb{Z}_K$ of $K$, it is well known that $\mathbb{Z}[\alpha]$ is a free $\mathbb{Z}$-module of rank $n$, from which it follows that the index $(\mathbb{Z}_K: \mathbb{Z}[\alpha])$ is finite. A well known formula linking $(\mathbb{Z}_K: \mathbb{Z}[\alpha])$, $\Delta(\alpha)$, and $d_K$ is given by:
	\begin{equation}\tag{1.1}
		\Delta(\alpha)=\pm (\mathbb{Z}_K: \mathbb{Z}[\alpha])^2\times d_K,
	\end{equation}  
	where $d_K$ is the absolute discriminant of $K$ and $\Delta(\alpha)$ is the discriminant of the minimal polynomial of $\alpha$ {over $\mathbb{Q}$}. The index of $K$, denoted by $i(K)$, is the greatest common divisor of the indices of all integral primitive elements of $K$. Say $i(K)=\gcd \ \{ ( \mathbb{Z}_K:\mathbb{Z}[\theta]) \, |\,{K=\mathbb{Q}(\theta)\text{ and }\theta \in \mathbb{Z}_K}\}$. A rational prime $p$ dividing $i(K)$ is called a prime common index divisor of $K$. {If $K$ is monogenic, then $\mathbb{Z}_K$ has a power integral basis, i.e., a $\mathbb{Z}$-basis of the form $(1,\theta,\dots,\theta^{n-1})$, and the index of $K$ is trivial, say $i(K)=1$. Therefore a field having a prime common index divisor is not monogenic.} In $1930$, Engstrom \cite{En} was the first one who studied the prime power decomposition of the index of a number field. He showed that for number fields of degree $n\leq 7$,  $\nu_p(i(K))$ is determined {by the form of the factorization of $p\mathbb{Z}_K$,  where $\nu_p$ is the $p$-adic valuation of $\mathbb{Q}$. It is an interesting problem to classify these number fields with non-trivial index, which are of course not monogenic}. In \cite{Sl}, \'Sliwa showed that, if $p$ is a non-ramified ideal in $K$, then $\nu_p(i(K))$ is determined by the factorization of $p\mathbb{Z}_K$. These results were generalized by Nart (\cite{N}), who developed a $p$-adic characterization of the index of a number field. In \cite{Nak}, Nakahara studied the index of non-cyclic but abelian biquadratic number fields. In \cite{GPP}  Ga\'al et al. characterized the field indices of biquadratic number fields having Galois group $V_4$. In \cite{DS}, for any quartic number field $K$ defined by a trinomial $x^4+ax+b$, Davis and Spearman characterized when $p=2,3$ divides $i(K)$. In \cite{EG}, for any  quartic number field $K$ defined by a trinomial $x^4+ax^2+b$, El Fadil and Ga\'al gave necessary and sufficient conditions on $a$ and $b$, which characterize when a rational prime $p$ divides $i(K)$. In \cite{E5}, for any rational prime $p$, El Fadil characterized when $p$ divides the index $i(K)$ for any quintic number field $K$ defined by a trinomial $x^5+ax^2+b$. In \cite{EK7}, for every rational prime $p$, we characterized $\nu_p(i(K))$ for any septic number field defined by a trinomial $x^7+ax^3+b$. In this paper, for any nonic number field $K$ defined by a monic irreducible trinomial $x^9+ax+b\in \mathbb{Z}[x]$ and for every rational prime $p$, we characterize when $p$ divides the index $i(K)$. Based on Engstrom's results given in \cite{En}, we evaluate $\nu_p(i(K))$ in some cases.  
	\section{Main Results}
	Throughout this section, $K$ is a number field generated by a  {complex} root $\alpha$ of a {monic irreducible} trinomial $F(x)=x^9+ax+b\in\mathbb{Z}[x]$. Without loss of generality, we assume that for every rational prime $p$,  $\nu_p(a)\leq 7$ or  $\nu_p(b) \leq 8$.\\
	\smallskip
	We start with the following theorem, which characterizes when the ring $\mathbb{Z}[\alpha]$ is integrally closed. 
	\begin{thm}\label{thm1}
		The ring $\mathbb{Z}[\alpha]$ is integrally closed if and only if the following conditions hold:
		\begin{enumerate}
			\item If $p$ divides both $a$ and $b$, then $\nu_p(b)=1$.
			\item If $2$ does not divide $a$ and divides $b$, then $(a,b)\in\{(1,0),(3,2)\}\md{4}$.
			\item If $3$ divides $a$ and does not divide $b$, then $$(a,b)\in \{(0,2),(0,5),(3,-1),(3,2)(6,-1),(6,5),(0,4),(0,7),(3,1),(3,7),(6,1),(6,4)\}\md{9}.$$
			\item For every rational prime $p\not\in\{2,3\}$ dividing $2^{24}a^9+3^{18}b^8$, if $\nu_p(ab)=0$, then  {$\nu_p(2^{24}a^9+3^{18}b^8)=1$}.
		\end{enumerate}
		{If all these conditions hold, then $K$ is monogenic and  $i(K)=1$}.
	\end{thm}
	
	\smallskip
	In the next theorems for every rational prime $p$, we characterize when $p$ divides the index $i(K)$ and  {evaluate $\nu_p(i(K))$ in some cases}.\\
	{Let us denote by $\Delta$ the discriminant of $F(x)$ and for every rational prime $p$, let $\Delta_p=\cfrac{\Delta}{p^{\nu_p(\Delta)}}$}.

	\smallskip 
	{The following theorem characterizes when $2$ divides $i(K)$}. 
	\begin{thm}\label{thmp2}
		The rational prime $2$ divides the index $i(K)$ if and only if one of the following conditions is satisfied:
		\begin{enumerate}
			\item $(a,b)\equiv(1,2)\md{4}$.
			
			\item $(a,b)\equiv(3,4)\md{8}$.
			
			\item $(a,b)\in\{(15,0),(7,8)\}\md{16}$.
			
			\item $(a,b)\equiv(28,0)\md{32}$.
			
			\item $(a,b)\in\{(4,0),(52,32)\}\md{64}$.
			
			\item $a\equiv112\md{128}$ and $b\equiv128\md{256}$.
			
			\item $(a,b)\in\{(368,256),(112,256),(240,0),(496,0),(448,0)\}\md{512}$.
			
			\item $a\equiv240\md{256}$ and $b\equiv256\md{512}$.
			
			\item $(a,b)\equiv(64,0)\md{1024}$.
		\end{enumerate}
		In particular, if one of the above conditions holds, then $K$ is not monogenic.
	\end{thm}
    \begin{rem}
    	Based on Engstrom's results given in \cite{En}, the following table provides $\nu_2(i(K))$ for some cases of Theorem $\ref{thmp2}$:
    	\begin{table}[htbp]
    		\centering
    		\caption{$\nu_2(i(K))$}
    		\begin{tabular}{|l|c|}
    			\hline
    			\multicolumn{1}{|c|}{\textbf{Conditions}}&\textbf{$\nu_2(i(K))$}\\
    			\hline
    			$(a,b)\equiv(1,2)\md{4}$&$1$\\
    			\hline
    			$(a,b)\equiv(7,8)\md{16}$ and $\nu_2(\Delta)$ is odd&\multirow{3}{*}{$3$}\\
    			\cline{1-1}
    			$(a,b)\equiv(7,8)\md{16}$, $\nu_2(\Delta)=28$, and $\Delta_2\equiv3\md{4}$&\\
    			\cline{1-1}
    			$(a,b)\equiv(7,8)\md{16}$, $\nu_2(\Delta)=2k\geq30$, and $\Delta_2\equiv 1\md{4}$&\\
    			\hline
    			$(a,b)\equiv(368,256)\md{512}$&$1$\\
    			\hline
    			$a\equiv240\md{256}$ and $b\equiv256\md{512}$ &$3$\\
    			\hline
    			\end{tabular}
    	\end{table}
    \end{rem}
	\smallskip
	The following theorem characterizes when $3$ divides $i(K)$.
	\begin{thm}\label{thmp3} 
	 The following table provides $\nu_3(i(K))$:
	 \begin{table}[htbp]
	 	\centering
	 	\caption{$\nu_3(i(K))$}
	 	\begin{tabular}{|l|l|c|}
	 		\hline
	 		\multicolumn{2}{|c|}{\textbf{Conditions}}&\textbf{$\nu_3(i(K))$}\\
	 		\hline
	 		$(a,b)\in \{(18,62),(72,8)\}\md{81}$ and $a+b\equiv -1\md{243}$&\multirow{2}{*}{$\nu_3(\Delta)=2k$ and}&\multirow{4}{*}{$1$}\\
	 		\cline{1-1}
	 		$(a,b)\equiv (45,35)\md{81}$ and $a+b\equiv 161\md{243}$&&\\
	 		\cline{1-1}
	 		$(a,b)\equiv \{(18,19),(72,73)\}\md{81}$ and $b-a\equiv 1\md{243}$&\multirow{2}{*}{$\Delta_3\equiv -1\md{3}$}&\\
	 		\cline{1-1}
	 		$(a,b)\equiv (45,46)\md{81}$ and $b-a\equiv 82\md{243}$&&\\
	 		\hline 
	 		\multicolumn{2}{|c|}{\textbf{Otherwise}}&$0$\\
	 		\hline
	 	\end{tabular}
	 \end{table}\\
	In particular, if $i(K)\neq1$, then $K$ is not monogenic.
	\end{thm}
	\smallskip 
	\begin{thm}\label{thmp}
		For every rational prime $p\geq 5$ and {for} every $(a,b)\in\mathbb{Z}^2$ such that $x^9+ax+b$ is irreducible, $p$ does not divide the index $i(K)$, where $K$ is the number field defined by $x^9+ax+b$.
	\end{thm}
	
	\section{Preliminaries}\label{sec3}
	{Our proofs are based on Newton polygon techniques applied on prime ideal factorization, which is a standard method which is rather technical but very efficient to apply. We have introduced the corresponding concepts in several former papers. Here we only give the theorem of index of Ore which plays a key role for proving our main results.}\\
	Let $K=\mathbb{Q}(\alpha)$ be a number field generated by a complex root $\alpha$ of a monic irreducible polynomial $F(x)\in\mathbb{Z}[x]$. We shall use Dedekind’s theorem \cite[Chapter I, Proposition 8.3]{Neu} and Dedekind’s criterion \cite[Theorem 6.1.4]{Co}. Let $\phi\in\mathbb{Z}_p[x]$ be a monic lift to an irreducible factor of $F(x)$ modulo $p$, $F(x)=a_0(x)+a_1(x)\phi(x)+\cdots+a_l(x)\phi(x)^l$ the $\phi$-expansion of $F(x)$, and $N_{\phi}^+(F)$ the principal $\phi$-Newton polygon of $F(x)$. Let $\mathbb{F}_{\phi}$ be the field $\mathbb{F}_p[x]/(\overline{\phi})$, then to every side $S$ of $N_{\phi}^+(F)$ with initial point $(i,u_i)$, and every $i=0,\ldots,l$, let the residue coefficient $c_i\in\mathbb{F}_{\phi}$ defined as follows: 
	$$c_{i}=
	\left
	\{\begin{array}{ll} 0,& \mbox{ if } (s+i,{\it u_{s+i}}) \mbox{ lies strictly
			above } S,\\
		\left(\dfrac{a_{s+i}(x)}{p^{{\it u_{s+i}}}}\right)
		\,\,
		\mod{(p,\phi(x))},&\mbox{ if }(s+i,{\it u_{s+i}}) \mbox{ lies on }S.
	\end{array}
	\right.$$
	Let $-\lambda=-h/e$ be the slope of $S$, where $h$ and $e$ are two positive coprime integers and $l=l(S)$ its length. Then  $d=l/e$ is the degree of $S$. Hence, if $i$ is not a multiple of $e$, then  $(s+i, u_{s+i})$ does not lie on $S$, and so $c_i=0$. Let ${R_{1}(F)(y)}=t_dy^d+t_{d-1}y^{d-1}+\cdots+t_{1}y+t_{0}\in\mathbb{F}_{\phi}[y]$, called  
	the residual polynomial of $F(x)$ associated to the side $S$, where for every $i=0,\dots,d$,  $t_i=c_{ie}$. If ${R_{1}(F)(y)}$ is square free for each side of the polygon $N_{\phi}^+(F)$, then we say that $F(x)$ is $\phi$-regular.\\ 
	Let $\overline{F(x)}=\prod_{i=1}^{r}\ol{\phi_i}^{l_i}$ be the factorization of $F(x)$ into powers of monic irreducible coprime polynomials over $\mathbb{F}_p$, we say that the polynomial $F(x)$ is $p$-{regular} if $F(x)$ is a $\phi_i$-regular polynomial with respect to $p$ for every $i=1,\dots,r$. Let  $N_{\phi_i}^+(F)=S_{i1}+\cdots+S_{ir_i}$ be the $\phi_i$-principal Newton polygon of $F(x)$ with respect to $p$. For every $j=1,\dots,r_i$, let ${R_{1_{ij}}(F)(y)}=\prod_{s=1}^{s_{ij}}\psi_{ijs}^{a_{ijs}}(y)$ be the factorization of {$R_{1_{ij}}(F)(y)$} in $\mathbb{F}_{\phi_i}[y]$. Then we have the following  theorem of index of Ore:
	\begin{thm}\label{thm4}$($\cite[Theorem 1.7 and Theorem 1.9]{EMN}$)$\\
		Under the above hypothesis, we have the following:
		\begin{enumerate}
			\item 
			$$\nu_p((\mathbb{Z}_K:\mathbb{Z}[\alpha]))\geq\sum_{i=1}^{r}\text{ind}_{\phi_i}(F).$$  
			The equality holds if $F(x) \text{ is }p$-regular.
			\item 
			If $F(x) \text{ is }p$-regular, then
			$$p\mathbb{Z}_K=\prod_{i=1}^r\prod_{j=1}^{t_i}\prod_{s=1}^{s_{ij}}\mathfrak{p}_{ijs}^{e_{ij}}$$
			is the factorization of $p\mathbb{Z}_K$ into powers of prime ideals of $\mathbb{Z}_K$, where $e_{ij}$ is the smallest positive integer satisfying $e_{ij}\la_{ij}\in \Z$ and the residue degree of $\mathfrak{p}_{ijs}$ over $p$ is given by $f_{ijs}=\deg(\phi_i)\times \deg(\psi_{ijs})$ for every $(i,j,s)$.
		\end{enumerate}
	\end{thm}
	\smallskip
	{If the theorem of Ore fails, that is, $F(x)$ is not $p$-regular, then in order to complete the factorization of $F(x)$, Guàrdia, Montes, and Nart introduced the notion of high order Newton polygon (\cite{GMN}). Similar to first order, for each order $r$, they introduced the valuation $\omega_r$, the key polynomial $\phi_r(x)$ for this valuation, the Newton polygon $N_r(F)$ of any polynomial $F(x)$ with respect to $\omega_r$ and $\phi_r(x)$, and for each side $T_i$ of $N_r(F)$, the residual polynomial $R_{r}(F)(y)$, and the index of $F(x)$ in order $r$. For more details, we refer to \cite{GMN}.}
	
	\section{Proofs of Main Results}
	Throughout this section, if $p\mathbb{Z}_K=\prod_{i=1}^r\prod_{j=1}^{t_i}\prod_{s=1}^{s_{ij}}\mathfrak{p}_{ijs}^{e_{ij}}$, then $e_{ij}$ denotes the ramification index of $\mathfrak{p}_{ijs}$ and $f_{ijs}$ denotes its residue degree for every $(i,j,s)$.\\
	For every rational prime $p$ and every integer $m$ let $m_p=m/p^{\nu_p(m)}$.\\
	
	\textit{\textbf{Proof of Theorem \ref{thm1}}}.\\
	Let $K=\mathbb{Q}(\alpha)$ be a number field defined by a monic irreducible trinomial $F(x)=x^{9}+ax+b\in\mathbb{Z}[x]$. Since  $\Delta=2^{24}a^9+3^{18}b^8$ is the discriminant of $F(x)$, thanks to the index formula $(1.1)$, we have the following:
	\begin{enumerate}
		\item If $p$ divides both $a$ and $b$, then $p$ does not divide the index $(\mathbb{Z}_K:\mathbb{Z}[\alpha])$ if and only if $\nu_p(b)=1$.
		
		\item If $p=2$ and $2$ does not divide $b$, then $2$ does not divide $(\mathbb{Z}_K:\mathbb{Z}[\alpha])$.
		
		\item If $p=2$, $2$ divides $b$ and does not divide $a$, then $F(x)\equiv x(x-1)^8\md{2}$. Let $\phi_1=x$ and $\phi_2=x-1$, then $F(x)=\cdots+36\phi_2^2+(a+9)\phi_2+a+b+1$. Hence $2$ does not divide $(\mathbb{Z}_K:\mathbb{Z}[\alpha])$ if and only if $\nu_2(a+b+1)=1$; that is $(a,b)\in\{(1,0),(3,2)\}\md{4}$.
		
		\item If $p=3$ and $3$ does not divide $a$, then $3$ does not divide $(\mathbb{Z}_K:\mathbb{Z}[\alpha])$.
		
		\item If $p=3$, $3$ divides $a$, and $b\equiv -1\md{3}$, then $F(x)\equiv (x-1)^9\md{3}$.  Let $\phi=x-1$. Then $F(x)=\cdots+36\phi^2+(a+9)\phi+a+b+1$. Hence $3$ does not divide $(\mathbb{Z}_K:\mathbb{Z}[\alpha])$ if and only if $\nu_3(a+b+1)=1$; that is $(a,b)\in\{(0,2),(0,5),(3,-1),(3,2)(6,-1),(6,5)\}\md{9}$.
		
		\item If $p=3$, $3$ divides $a$, and $b\equiv 1\md{3}$, then $F(x)\equiv (x+1)^9\md{3}$. Let $\phi=x+1$. Then
		$F(x)=\cdots-36\phi^2+(a+9)\phi-a+b-1$. Hence $3$ does not divide $(\mathbb{Z}_K:\mathbb{Z}[\alpha])$ if and only if $\nu_3(-a+b-1)=1$; that is $(a,b)\in\{(0,4),(0,7),(3,1),(3,7),(6,1),(6,4)\}\md{9}$.
		
		\item If $p\notin\{2,3\}$, $p^2$ divides $\Delta=2^{24}a^9+3^{18}b^8$, and $\nu_p(ab)=0$, then $F(x)$ admits a multiple root $\overline{u}$ modulo $p$ {if and only if} $F(u)=u^{9}+au+b\equiv 0\md{p}$ and $F'(u)=9u^{8}+a\equiv 0\md{p}$. That is {$8au+9b\equiv 0\md{p}$}. Let $u=\cfrac{-9b}{8a}\in\mathbb{Q}$. Since $\nu_p(8a)=0$, then $u\in\mathbb{Z}_p$. Let $\phi=x-u$. Then $F(x)=\cdots+36u^{7}\phi^2+A\phi+B$, with
		$$
		\begin{array}{lllllll}
			A&=&a+9u^8&=&\cfrac{2^{24}a^7+3^{18}b^8}{2^{24}a^8}&=&\cfrac{\Delta}{2^{24}a^8}~~, {\text{and}}\\
			B&=&au+b+u^9&=&\cfrac{-3^{18}b^9-2^{24}a^9b}{2^{27}a^9}&=&\cfrac{-b\Delta}{2^{27}a^9}~~.
		\end{array}
		$$ 
		Since $\nu_p(A)=\nu_p(B)=\nu_p(\Delta)$ and $\nu_p(36u^7)=0$, then $N_{\phi}^+(F)=S_1$ has a single side joining $(0,\nu_p(\Delta))$ and $(2,0)$. Since $\nu_p(\Delta)\geq2$, by Theorem $\ref{thm4}$, $\nu_p((\mathbb{Z}_K:\mathbb{Z}[\alpha]))\geq \lfloor \nu_p(\Delta)/2\rfloor\geq1$. Hence $p$ divides the index $(\mathbb{Z}_K:\mathbb{Z}[\alpha])$.		
	\end{enumerate}
	\begin{flushright}
		$\square$
	\end{flushright}
	
	For the proof of Theorems \ref{thmp2}, \ref{thmp3}, and \ref{thmp}, we need the following lemma, which characterizes the prime common index divisors of $K$.
	\begin{lem}\label{index}$($\cite{En}$)$\\
		Let $p$ be a rational prime and {$K$ a} number field. For every positive integer $f$, let $\mathcal{P}_f$ be the number of distinct prime ideals of $\mathbb{Z}_K$ lying above $p$ with residue degree $f$ and $\mathcal{N}_f$ the number of monic irreducible polynomials of $\mathbb{F}_p[x]$ of degree $f$.  Then $p$ {divides the index $i(K)$} if and only if $\mathcal{P}_f>\mathcal{N}_f$ for some positive integer $f$.
	\end{lem}
	\smallskip	
	
	\textit{\textbf{Proof of Theorem \ref{thmp2}}}.\\
	Since $\Delta=2^{24}a^9+3^{18}b^8$ is the discriminant of $F(x)$, thanks to the index formula $(1.1)$, if $2$ does not divide $b$, then $\nu_2((\mathbb{Z}_K:\mathbb{Z}[\alpha]))=0$ and so $\nu_2(i(K))=0$. Assume that $2$ divides $b$. {Then} we have the following cases:
	\begin{enumerate}
		\item If $\nu_2(a)=0$, then $F(x)\equiv x(x-1)^8\md{2}$. Let {$\phi=x-1$. Then $x$} provides a unique prime ideal of $\mathbb{Z}_K$ lying above $2$ with residue degree $1$. For {$\phi$, let} {$F(x)=\phi^9+9\phi^8+36\phi^7+84\phi^6+126\phi^5+126\phi^4+84\phi^3+36\phi^2+(a+9)\phi+a+b+1$. Then} we have the following cases:
		\begin{enumerate}
			\item[(i)] If $(a,b)\in\{(1,0),(3,2)\}\md{4}$, then by Theorem $\ref{thm1}$, $\nu_2(i(K))=0$.
			
			\item[(ii)] If $(a,b)\equiv (1,2)\md{4}$, then $N_{\phi}^+(F)=S_{1}+S_{2}$ has two sides joining $(0,w)$, $(1,1)$, and $(8,0)$ with $w\geq2$. Thus the degree of each side is $1$ and so $2\mathbb{Z}_K=\mathfrak{p}_{11}\mathfrak{p}_{21}\mathfrak{p}_{22}^7$ with residue degree $1$ each {ideal factor}. Since there are just two monic irreducible polynomials of degree $1$ in $\mathbb{F}_2[x]$,  by Lemma $\ref{index}$, $2$ divides $i(K)$. Applying \cite[Corollary, p: 230]{En}, we get $\nu_2(i(K))=1$.
			
			\item[(iii)] For $(a,b)\equiv (3,0)\md{4}$, we have the following sub-cases:
			\begin{enumerate}
				\item[(a)] {If $(a,b)\not\equiv (7,8)\md{16}$, then the treatment of this case is similar the case (ii) above, and Table $3$ summarizes the obtained results.}
				\begin{table}[htbp]
					\centering
					\caption{}
					\begin{tabular}{|l|c|l|c|}
						\hline
						\multicolumn{1}{|c|}{Cases}&$2\mathbb{Z}_K$&\multicolumn{1}{|c|}{$f_i$}&$\nu_2(i(K))$\\
						\hline
						$(a,b)\in\{(3,0),(7,4)\}\md{8}$&$\mathfrak{p}_1\mathfrak{p}_2^4$&$f_1=1$, $f_2=2$&$0$\\
						\hline
						$(a,b)\equiv (3,4)\md{8}$&$\mathfrak{p}_{1}\mathfrak{p}_{2}\mathfrak{p}_{3}^3\mathfrak{p}_{4}^4$&$f_i=1$&$\geq1$\\
						\hline
						$(a,b)\in\{(7,0),(15,8)\}\md{16}$&$\mathfrak{p}_{1}\mathfrak{p}_{2}^2\mathfrak{p}_{3}^4$&$f_{1}=f_{3}=1$, $f_{2}=2$&$0$\\
						\hline
						$(a,b)\equiv (15,0)\md{16}$&$\mathfrak{p}_{1}\mathfrak{p}_{2}\mathfrak{p}_{3}^3\mathfrak{p}_{4}^4$&$f_i=1$&$\geq1$\\
						\hline
					\end{tabular}
				\end{table}
				
				\item[(b)] If $(a,b)\equiv (7,8)\md{16}$, then $\nu_2(\Delta)\geq28$. As in the proof of Theorem $\ref{thm1}$, let $b_2=\cfrac{b}{8}$ and $u=\cfrac{-9b_2}{a}$. Since $2$ does not divide $a$, then $u\in\mathbb{Z}_2$. Let $\phi=x-u$, then $F(x)=\phi^9+9u\phi^8+36u^2\phi^7+84u^3\phi^6+126u^4\phi^5+126u^5\phi^4+84u^6\phi^3+36u^7\phi^2+A\phi+B$ with
				$$
				\begin{array}{lllllll}
					A&=&a+9u^8&=&\cfrac{\Delta}{2^{24}a^8},\text{ and}&&\\
					B&=&au+b+u^9&=&\cfrac{-b\Delta}{2^{27}a^9}&=&\cfrac{-b_2\Delta}{2^{24}a^9}.
				\end{array}
				$$
				Thus $\nu_2(A)=\nu_2(B)=\nu_2(\Delta)-24$. 
				\item[b$_1$-] If $\nu_2(\Delta)$ is odd, then $N_{\phi}^+(F)=S_{1}+S_{2}+S_{3}$ has three sides joining $(0,\nu_2(\Delta)-24)$, $(2,2)$, $(4,1)$, and $(8,0)$. Thus the degree of each side is $1$ and so $2\mathbb{Z}_K=\mathfrak{p}_{11}\mathfrak{p}_{21}^2\mathfrak{p}_{22}^2\mathfrak{p}_{23}^4$ with residue degree $1$ each {ideal factor}. Hence $2$ divides $i(K)$. Applying \cite[Theorem 6]{En}, we get $\nu_2(i(K))=3$. 
				\item[b$_2$-] If $\nu_2(\Delta)=2k+26$ is even $(k\geq1)$, then $N_{\phi_2}^+(F)=S_{1}+S_{2}+S_{3}$ has three sides joining $(0,2k+2)$ $(2,2)$, $(4,1)$, and $(8,0)$ with $d(S_2)=d(S_3)=1$ and ${R_{1_1}(F)(y)}=(y+1)^2\in\mathbb{F}_{\phi}[y]$. Let us replace $\phi$ by $x-(u+2^k)$. Then $F(x)=\cdots+A_4(x-(u+2^k))^4+A_3(x-(u+2^k))^3+A_2(x-(u+2^k))^2+A_1(x-(u+2^k))+A_0$ with  $\nu_2(A_4)=1$, $\nu_2(A_3)=\nu_2(A_2)=2$, 
				$$
				\begin{array}{lll}
					A_1&=&A+72u^72^k+252u^6(2^k)^2+504u^5(2^k)^3+630u^4(2^k)^4+504u^3(2^k)^5+252u^2(2^k)^6\\
					&&+72u(2^k)^7+9(2^k)^8, {\text{ and}}\\
					A_0&=&B+2^kA+36u^7(2^k)^2+84u^6(2^k)^3+126u^5(2^k)^4+126u^4(2^k)^5\\
					&&+84u^3(2^k)^6+36u^2(2^k)^7+9u(2^k)^8+(2^k)^9.
				\end{array}
				$$
				{Clearly, $\nu_2(A_1)=k+3$}.\\
				\item[-] For $\nu_2(\Delta)=28$; $k=1$, we have $A_0\equiv B+2A+36u^72^2+84u^62^3+126u^52^4+126u^42^5\md{128}$. Hence $A_0\equiv\cfrac{2^4}{a^9}(-b_2\Delta_2+2a\Delta_2+63a^2b_2^7+10a^3b_2^6+82a^4b_2^5+124a^5b_2^4)\md{128}$. {Since $\cfrac{A_0}{2^4}\equiv \cfrac{1}{a^9}(-b_2\Delta_2+2a\Delta_2+7a^2b_2^7+2a^3b_2^6+2a^4b_2^5+4a^5b_2^4)\md{8}\equiv-(-b_2\Delta_2+6\Delta_2+b_2+2)\md{8}$, three cases arise are summarized in Table $4$.}
				\begin{table}[htbp]
					\centering
					\caption{$\nu_2(\Delta)=28$}
					\begin{tabular}{|l|c|c|l|c|}
						\hline
						\multicolumn{1}{|c|}{Cases}&$\nu_2(A_0)$&$2\mathbb{Z}_K$&\multicolumn{1}{|c|}{$f_i$}&$\nu_2(i(K))$\\
						\hline
						$\Delta_2\equiv 1\md{8}$&$\geq7$&$\mathfrak{p}_{1}\mathfrak{p}_{2}\mathfrak{p}_{3}\mathfrak{p}_{4}^2\mathfrak{p}_{5}^4$&$f_i=1$&\multirow{2}{*}{$\geq1$}\\
						\cline{1-4} 
						$\Delta_2\equiv 5\md{8}$&$6$&$\mathfrak{p}_{1}\mathfrak{p}_{2}\mathfrak{p}_{3}^2\mathfrak{p}_{4}^4$&$f_{1}=f_{3}=f_{4}=1$, $f_{2}=2$&\\
						\hline
						$\Delta_2\equiv 3,7\md{8}$&$5$&$\mathfrak{p}_{1}\mathfrak{p}_{2}^2\mathfrak{p}_{3}^2\mathfrak{p}_{4}^4$&$f_i=1$&$3$\\
						\hline 
					\end{tabular}
				\end{table}
				\item[-] For $\nu_2(\Delta)\geq30$; $k\geq2$, we have $A_0\equiv B+36u^7(2^k)^2\md{2^{3k+3}}$. Hence $A_0\equiv\cfrac{2^{2k+2}}{a^9}(-b_2\Delta_2+9^8a^3b_2^7))\md{2^{3k+3}}$. {Since $\cfrac{A_0}{2^{2k+2}}\equiv -(-b_2\Delta_2+7b_2)\md{8}$, three cases arise are summarized in Table $5$.}
				\begin{table} [htbp]
					\centering
					\caption{$\nu_2(\Delta)\geq30$ is even}
					\begin{tabular}{|l|c|c|l|c|}
						\hline
						\multicolumn{1}{|c|}{Cases}&$\nu_2(A_0)$&$2\mathbb{Z}_K$&\multicolumn{1}{|c|}{$f_i$}&$\nu_2(i(K))$\\
						\hline
						$\Delta_2\equiv 7\md{8}$,&$\geq2k+5$&$\mathfrak{p}_{1}\mathfrak{p}_{2}\mathfrak{p}_{3}\mathfrak{p}_{4}^2\mathfrak{p}_{5}^4$&$f_i=1$&\multirow{2}{*}{$\geq1$}\\
						\cline{1-4}  
						$\Delta_2\equiv 3\md{8}$&$2k+4$&$\mathfrak{p}_{1}\mathfrak{p}_{2}\mathfrak{p}_{3}^2\mathfrak{p}_{4}^4$&$f_{1}=f_{3}=f_{4}=1$, $f_{2}=2$&\\
						\hline
						$\Delta_2\equiv 1,5\md{8}$&$2k+3$&$\mathfrak{p}_{1}\mathfrak{p}_{2}^2\mathfrak{p}_{3}^2\mathfrak{p}_{4}^4$&$f_i=1$&$3$\\
						\hline 
					\end{tabular}
				\end{table}
				
			\end{enumerate}
		\end{enumerate}
		\item If $\nu_2(a)\geq1$, then $F(x)\equiv x^9\md{2}$. Let $\phi=x$. Then $F(x)=\phi^9+a\phi+b$. By assumption, $\nu_2(b)\leq8$ or $\nu_2(a)\leq7$.\\
		If $8\nu_2(b)<9\nu_2(a)$, then $N_{\phi}(F)=S_1$ has a single side of degree $d\in\{1,3\}$.
		\begin{enumerate}
			\item[(i)] If $d=1$ then {$R_{1_1}(F)(y)$} is irreducible as it is of degree $1$. Thus $2\mathbb{Z}_K=\mathfrak{p}_1^{9}$ with residue degree $1$. Hence $\nu_2(i(K))=0$.
			
			\item[(ii)] If $d=3$, then ${R_{1_1}(F)(y)}=y^3+1=(y+1)(y^2+y+1)\in\mathbb{F}_{\phi}[y]$. Thus $2\mathbb{Z}_K=\mathfrak{p}_1^3\mathfrak{p}_2^3$ with $f_1=1$ and $f_2=2$. Hence $\nu_2(i(K))=0$.
		\end{enumerate}
		
		If $8\nu_2(b)>9\nu_2(a)$, then {$N_{\phi}(F)=S_1+S_2$ has two sides joining $(0,\nu_2(b))$, $(1,\nu_2(a))$, and $(0,9)$ with $d(S_1)=1$ and $d(S_2)\in\{1,2,4\}$ since $\nu_2(a)\leq 7$. Let $d=d(S_2)$. Then we have the following cases:}
		\begin{enumerate}
			\item[(i)] If $d=1$; $\nu_2(a)\in\{1,3,5,7\}$, then $2\mathbb{Z}_K=\mathfrak{p}_1\mathfrak{p}_2^8$ with residue degree $1$ each {ideal factor}. Hence $\nu_2(i(K))=0$.
			\item[(ii)] If $d=2$ with $\nu_2(b)\geq3$ and $\nu_2(a)=2$; $a\equiv 4\md{8}$ and $b\equiv 0\md{8}$, then $N_{\phi}(F)=S_1+S_2$ has two sides joining $(0,\nu_2(b))$, $(1,2)$ and $(9,0)$ with $d(S_1)=1$ and ${R_{1_2}}(F)(y)=(y+1)^2\in\mathbb{F}_{\phi}[y]$. In this case, we have to use second order Newton polygon. Let $\omega_2=4[\nu_2,\phi,1/4]$ be the valuation of second order Newton polygon and $g_2=x^4-2$ the key polynomial of $\omega_2$, where $[\nu_2,\phi,1/4]$ is the augmented valuation of $\nu_2$ with respect to $\phi$ and $\lambda=1/4$. Let $F(x)=xg_2^2+4xg_2+(a+4)x+b$, then we have $\omega_2(x)=1$, $\omega_2(g_2)=4$, and $\omega_2(m)=4\times \nu_2(m)$ for every $m\in\mathbb{Q}_2$. {Table $6$ summarizes the obtained results.}
			
			\begin{table}[htbp]
				\centering
				\caption{}
				\begin{tabular}{|l|c|c|l|c|}
					\hline
					\multicolumn{1}{|c|}{Cases}&$g_2$&$2\mathbb{Z}_K$&\multicolumn{1}{|c|}{$f_i$}&$\nu_2(i(K))$\\
					\hline
					$(a,b)\in\{(4,8),(12,8)\}\md{16}$&\multirow{4}{*}{$x^4-2$}&\multirow{2}{*}{$\mathfrak{p}_1\mathfrak{p}_2^8$}&\multirow{2}{*}{$f_i=1$}&\multirow{3}{*}{$0$}\\
					\cline{1-1}
					$(a,b)\in\{(12,16),(28,16)\}\md{32}$&&&&\\
					\cline{1-1}\cline{3-4}
					$(a,b)\equiv (12,0)\md{32}$&&$\mathfrak{p}_1\mathfrak{p}_2^4$&$f_1=1$, $f_2=2$&\\
					\cline{1-1}\cline{3-5}
					$(a,b)\equiv (28,0)\md{32}$&&$\mathfrak{p}_{11}\mathfrak{p}_{21}^4\mathfrak{p}_{22}^4$&$f_i=1$&$\geq1$\\
					\hline
					$(a,b)\in\{(4,16),(20,16)\}\md{32}$&$x^4-2x^2-2$&\multirow{2}{*}{$\mathfrak{p}_1\mathfrak{p}_2^8$}&\multirow{2}{*}{$f_i=1$}&\multirow{3}{*}{$0$}\\
					\cline{1-2}
					
					$(a,b)\in\{(4,32),(36,32)\}\md{64}$&\multirow{3}{*}{$x^4-2x^2-6$}&&&\\
					\cline{1-1}\cline{3-4}
					$(a,b)\equiv (36,0)\md{64}$&&$\mathfrak{p}_1\mathfrak{p}_2^4$& $f_1=1$, $f_2=2$&\\
					\cline{1-1}\cline{3-5}
					$(a,b)\equiv (4,0)\md{64}$&&$\mathfrak{p}_{11}\mathfrak{p}_{21}^4\mathfrak{p}_{22}^4$&$f_i=1$&$\geq1$\\
					\hline
					$(a,b)\in\{(20,0),(52,0)\}\md{64}$& \multirow{3}{*}{$x^4-2x^2-4x-2$}&$\mathfrak{p}_1\mathfrak{p}_2^8$&$f_i=1$&\multirow{2}{*}{$0$}\\
					\cline{1-1}\cline{3-4}
					$(a,b)\equiv (20,32)\md{64}$&&$\mathfrak{p}_1\mathfrak{p}_2^4$& $f_1=1$, $f_2=2$&\\
					\cline{1-1}\cline{3-5}
					$(a,b)\equiv (52,32)\md{64}$&&$\mathfrak{p}_{11}\mathfrak{p}_{21}^4\mathfrak{p}_{22}^4$&$f_i=1$&$\geq1$\\
					\hline
				\end{tabular}
			\end{table}

			\item[(iii)] If $d=4$; $a\equiv 16\md{32}$ and $b\equiv 0\md{32}$, then $N_{\phi}(F)=S_1+S_2$ has two sides joining $(0,\nu_2(b))$, $(1,4)$ and $(9,0)$ with $d(S_1)=1$ and ${R_{1_2}}(F)(y)=(y+1)^4\in\mathbb{F}_{\phi}[y]$. In this case, we have to use second order Newton polygon. Let $\omega_2=2[\nu_2,\phi,1/2]$ be the valuation of second order Newton polygon and $g_2=x^2-2$ the key polynomial of $\omega_2$, where $[\nu_2,\phi,1/2]$ is the augmented valuation of $\nu_2$ with respect to $\phi$ and $\lambda=1/2$. Let $F(x)=xg_2^4+8xg_2^3+24xg_2^2+32xg_2+(a+16)x+b$, then we have $\omega_2(x)=1$, $\omega_2(g_2)=2$, and $\omega_2(m)=2\times \nu_2(m)$ for every $m\in\mathbb{Q}_2$.
			
			\begin{enumerate}
				\item[(a)] If $(a,b)\in\{(16,32),(48,32)\}\md{64}$, then $N_2^+(F)=T_1$ has a single side joining $(0,10)$ and $(4,9)$. Thus $2\mathbb{Z}_K=\mathfrak{p}_1\mathfrak{p}_2^8$ with residue degree $1$ each {ideal factor}. Hence $\nu_2(i(K))=0$. 
				
				\item[(b)] If $(a,b)\equiv (16,0)\md{64}$, then $N_2^+(F)=T_1$ has a single side joining $(0,11)$ and $(2,9)$ with ${R_{2_1}(F)(y)}=(y+1)^2\in\mathbb{F}_2[y]$. In this case, we have to use third order Newton polygon. Let $\omega_3=2[\omega_2,g_2,1/2]$ be the valuation of third order Newton polygon and $g_3=x^4-4x^2-4x+4$ the key polynomial of $\omega_3$, where $[\omega_2,g_2,1/2]$ is the augmented valuation of $\omega_2$ with respect to $g_2$, and $\lambda'=1/2$. Let 
				$$F(x)=xg_3^2+((8x+8)g_2+24x+48)g_3+(48x+160)g_2+(a+176)x+b+192,$$
				then we have $\omega_3(x)=2$, $\omega_3(g_2)=5$, $\omega_3(g_3)=10$, and $\omega_3(m)=4\times \nu_2(m)$ for every $m\in\mathbb{Q}_2$. Thus $N_3^+(F)=T'_1$ has a single side joining $(0,23)$ and $(2,22)$. It follows that $2\mathbb{Z}_K=\mathfrak{p}_1\mathfrak{p}_2^8$ with residue degree $1$ each {ideal factor}. Hence $\nu_2(i(K))=0$.
				
				\item[(c)] {For the other cases, we need also to use third order Newton polygon, and its treatment is similar to the case (b) above. The obtained results are summarized in Table $7$}.
				
				\begin{table}[htbp]
					\centering
					\caption{}
					\begin{tabular}{|l|c|c|c|l|c|}
						\hline
						\multicolumn{1}{|c|}{Cases}&$g_2$&$g_3$&$2\mathbb{Z}_K$&\multicolumn{1}{|c|}{$f_i$}&$\nu_2(i(K))$\\
						\hline
						$(a,b)\in\{(48,64),(112,64)\}$&\multirow{16}{*}{\rotatebox{90}{$x^2-2$\,}}&\multirow{2}{*}{$-$}&\multirow{2}{*}{$\mathfrak{p}_1\mathfrak{p}_2^8$}&\multirow{2}{*}{$f_i=1$}&\multirow{3}{*}{$0$}\\
						\multicolumn{1}{|c|}{$\md{128}$}&&&&&\\
						\cline{1-1}\cline{3-5}
						$(a,b)\equiv(48,0)\md{128}$&&$x^4-2x^3-4x^2+4x-4$&{$\mathfrak{p}_1\mathfrak{p}_2^4$}&{$f_1=1$, $f_2=2$}&\\
						\cline{1-1}\cline{3-6} 
						$(a,b)\in\{(112,128),(240,128)\}$&&\multirow{2}{*}{$x^2-2x-2$}&\multirow{2}{*}{$\mathfrak{p}_1\mathfrak{p}_2^4\mathfrak{p}_3^4$}&\multirow{2}{*}{$f_i=1$}&\multirow{2}{*}{$\geq1$}\\
						\multicolumn{1}{|c|}{$\md{256}$}&&&&&\\
						\cline{1-1}\cline{3-6}
						$(a,b)\in\{(112,0),(368,0)\}$&&\multirow{2}{*}{$x^2-2x-2$}&\multirow{2}{*}{$\mathfrak{p}_1\mathfrak{p}_2^2\mathfrak{p}_3^4$}&$f_1=f_3=1$&\multirow{2}{*}{$0$}\\
						\multicolumn{1}{|c|}{$\md{512}$}&&&&$f_2=2$&\\
						\cline{1-1}\cline{3-6}
						\multirow{2}{*}{$(a,b)\equiv (368,256)\md{512}$}&&\multirow{2}{*}{$x^2-2x-2$}&\multirow{2}{*}{$\mathfrak{p}_1\mathfrak{p}_2^2\mathfrak{p}_3^2$}& $f_1=1$&\multirow{2}{*}{$1$}\\
						&&&&$f_2=f_3=2$&\\
						\cline{1-1}\cline{3-6} 
						\multirow{2}{*}{$(a,b)\equiv (112,256)\md{512}$}&&\multirow{2}{*}{$x^2-2x-2$}&\multirow{2}{*}{$\mathfrak{p}_{1}\mathfrak{p}_{2}^2\mathfrak{p}_{3}^2\mathfrak{p}_{4}^2$}& $f_{1}=f_{3}=f_{4}=1$&\multirow{2}{*}{$\geq1$}\\
						&&&&$f_{2}=2$&\\
						\cline{1-1}\cline{3-6} 
						$(a,b)\in\{(240,256),(496,256)\}$&&\multirow{5}{*}{$x^2-2x-6$}&\multirow{2}{*}{$\mathfrak{p}_{1}\mathfrak{p}_{2}^2\mathfrak{p}_{3}^2\mathfrak{p}_{4}^4$}&\multirow{2}{*}{$f_i=1$}&\multirow{2}{*}{$3$}\\
						\multicolumn{1}{|c|}{$\md{512}$}&&&&&\\
						\cline{1-1}\cline{4-6} 
						
						\multirow{2}{*}{$(a,b)\equiv(240,0)\md{512}$}&&&\multirow{2}{*}{$\mathfrak{p}_1\mathfrak{p}_{2}^2\mathfrak{p}_{3}^2\mathfrak{p}_4^2$}& $f_1=f_{2}=f_{3}=1$&\multirow{2}{*}{$\geq1$}\\
						&&&&$f_4=2$&\\
						\cline{1-1}\cline{4-6} 
						$(a,b)\equiv (496,0)\md{512}$&&&$\mathfrak{p}_{1}\mathfrak{p}_{2}^2\mathfrak{p}_{3}^2\mathfrak{p}_{4}^2\mathfrak{p}_{5}^2$&$f_i=1$&$\geq1$\\
						\hline
					\end{tabular}
				\end{table}
				
			\end{enumerate}
			
			\item[(iv)] If $d=2$ with $\nu_2(b)\geq7$ and $\nu_2(a)=6$; $a\equiv 64\md{128}$ and $b\equiv 0\md{128}$, then $N_{\phi}^+(F)=S_1+S_2$ has two sides joining $(0,\nu_2(b))$, $(1,6)$, and $(9,0)$ with $d(S_1)=1$ and ${R_{1_2}(F)(y)}=(y+1)^2\in\mathbb{F}_{\phi}[y]$. In this case, we have to use second order Newton polygon. Let $\omega_2=4[\nu_2,\phi,3/4]$ be the valuation of second order Newton polygon and $g_2=x^4-8$ the key polynomial of $\omega_2$, where $[\nu_2,\phi,3/4]$ is the augmented valuation of $\nu_2$ with respect to $\phi$ and $\lambda=3/4$. Let $F(x)=xg_2^2+16xg_2+(a+64)x+b$, then we have $\omega_2(x)=3$, $\omega_2(g_2)=12$, and $\omega_2(m)=4\times \nu_2(m)$ for every $m\in\mathbb{Q}_2$. {Table $8$ summarizes the obtained results.}
			
			\begin{table}[htbp]
				\centering
				\caption{}
				\begin{tabular}{|l|c|c|l|c|}
					\hline
					\multicolumn{1}{|c|}{Cases}&$g_2$&$2\mathbb{Z}_K$&\multicolumn{1}{|c|}{$f_i$}&$\nu_2(i(K))$\\
					\hline
					$(a,b)\in\{(64,128),(192,128)\}\md{256}$&$x^4-8$&\multirow{3}{*}{$\mathfrak{p}_1\mathfrak{p}_2^8$}&\multirow{3}{*}{$f_i=1$}&\multirow{4}{*}{$0$}\\
					\cline{1-2}
					$(a,b)\in\{(64,256),(320,256)\}\md{512}$&$x^4-4x^2-8$&&&\\
					\cline{1-2} 
					$(a,b)\in\{(64,512),(566,512)\}\md{1024}$&$x^4-4x-24$&&&\\
					\cline{1-4}
					$(a,b)\equiv (576,0)\md{1024}$&\multirow{2}{*}{$x^4-4x-24$}&$\mathfrak{p}_1\mathfrak{p}_2^4$& $f_1=1$, $f_2=2$&\\
					\cline{1-1}\cline{3-5}
					$(a,b)\equiv (64,0)\md{1024}$&&$\mathfrak{p}_{1}\mathfrak{p}_{2}^4\mathfrak{p}_{3}^4$&$f_i=1$&$\geq1$\\
					\hline 
					$(a,b)\in\{(320,512),(832,512)\}\md{1024}$&$x^4-4x^2-8$&\multirow{3}{*}{$\mathfrak{p}_1\mathfrak{p}_2^8$}&\multirow{3}{*}{$f_i=1$}&\multirow{4}{*}{$0$}\\
					\cline{1-2}
					$(a,b)\in\{(320,0),(832,0)\}\md{1024}$&$x^4-2x^3-4x^2-8$&&&\\
					\cline{1-2}
					$(a,b)\in\{(192,256),(448,256)\}\md{512}$&\multirow{3}{*}{$x^4-8$}&&&\\
					\cline{1-1}\cline{3-4}
					$(a,b)\equiv (192,0)\md{512}$&&$\mathfrak{p}_1\mathfrak{p}_2^4$&$f_1=1$, $f_2=2$&\\
					\cline{1-1}\cline{3-5}
					$(a,b)\equiv (448,0)\md{512}$&&$\mathfrak{p}_{1}\mathfrak{p}_{2}^4\mathfrak{p}_{3}^4$&$f_i=1$&$\geq1$\\
					\hline
				\end{tabular}
			\end{table}
			
		\end{enumerate}
		
	\end{enumerate}
	
	\begin{flushright}
		$\square$
	\end{flushright}
	\textit{\textbf{Proof of Theorem \ref{thmp3}}}.\\
	{If $\nu_3(a)=0$, then since $\Delta=2^{24}a^9+3^{18}b^8$ is the discriminant of $F(x)$, thanks to the index formula $(1.1)$, $\nu_3((\mathbb{Z}_K:\mathbb{Z}[\alpha]))=0$ and so $\nu_3(i(K))=0$. Now assume that $3$ divides $a$. Then} we have the following cases:
	\begin{enumerate}
		\item If $b\equiv -1\md{3}$, then $F(x)\equiv (x-1)^9\md{3}$. Let $\phi=x-1$. Then $F(x)=\phi^9+9\phi^8+36\phi^7+84\phi^6+126\phi^5+126\phi^4+84\phi^3+36\phi^2+(a+9)\phi+a+b+1$.
		\begin{enumerate}
			\item[(i)] If $(a,b)\in\{(0,2),(0,5),(3,-1),(3,2)(6,-1),(6,5)\}\md{9}$, then by Theorem $\ref{thm1}$, $\nu_3(i(K))=0$.
			
			\item[(ii)] If $(a,b)\in\{(3,5),(6,2)\}\md{9}$, then $N_{\phi}(F)=S_1+S_2$ has two sides joining $(0,w)$, $(1,1)$, and $(9,0)$ with $w\geq2$. Thus the degree of each side is $1$ {and so} $3\mathbb{Z}_K=\mathfrak{p}_1\mathfrak{p}_2^8$ with residue degree $1$ each {ideal factor}. Hence $\nu_3(i(K))=0$.
			
			\item[(iii)] If $(a,b)\in\{(0,8),(0,17),(9,-1),(9,8),(18,-1),(18,17)\}\md{27}$, then  $N_{\phi}(F)=S_1+S_2$ has two sides joining $(0,2)$, $(3,1)$, and $(9,0)$. Thus $3\mathbb{Z}_K=\mathfrak{p}_1^3\mathfrak{p}_2^6$ with residue degree $1$ each {ideal factor}. Hence $\nu_3(i(K))=0$.
			
			\item[(iv)] If $(a,b)\in\{(0,-1),(9,17)\}\md{27}$, then $N_{\phi}(F)=S_1+S_2+S_3$ has three sides joining $(0,w)$, $(1,2)$, $(3,1)$, and $(9,0)$ with $w\geq3$. It follows that $3\mathbb{Z}_K=\mathfrak{p}_1\mathfrak{p}_2^2\mathfrak{p}_3^6$ with residue degree $1$ each {ideal factor}. Hence $\nu_3(i(K))=0$.
			
			\item[(v)] If $(a,b)\in\{(18,8),(18,35),(45,8),(45,62),(72,35),(72,62)\}\md{81}$, then $N_{\phi}(F)=S_1+S_2$ has two sides joining $(0,3)$, $(3,1)$, and $(9,0)$. It follows that $3\mathbb{Z}_K=\mathfrak{p}_1^3\mathfrak{p}_2^6$ with residue degree $1$ each {ideal factor}. Hence $\nu_3(i(K))=0$.
			
			\item [(vi)] If $(a,b)\equiv (18,62)\md{81}$, then we have the following sub-cases:
			\begin{enumerate}
				\item[(a)] If {$a+b\equiv 80\md{243}$}, then $N_{\phi}(F)=S_1+S_2$ has two sides joining $(0,4)$, $(3,1)$, and $(9,0)$ with $d(S_2)=1$ and ${R_{1_1}(F)(y)}=y^3+y^2+y+1=(y^2+1)(y+1)\in\mathbb{F}_{\phi}[y]$. Thus $3\mathbb{Z}_K=\mathfrak{p}_{11}\mathfrak{p}_{12}\mathfrak{p}_{21}^6$ with $f_{11}=f_{21}=1$ and $f_{12}=2$. Hence $\nu_3(i(K))=0$. 
				
				\item[(b)] If {$a+b\equiv 161\md{243}$}, then $N_{\phi}(F)=S_1+S_2$ has two sides joining $(0,4)$, $(3,1)$, and $(9,0)$ with $d(S_2)=1$ and ${R_{1_1}(F)(y)}=y^3+y^2+y-1$ which is irreducible over $\mathbb{F}_{\phi}$. Thus $3\mathbb{Z}_K=\mathfrak{p}_1\mathfrak{p}_2^6$ with $f_{1}=3$ and $f_{2}=1$. Hence $\nu_3(i(K))=0$.
				
				\item[(c)] If {$a+b\equiv -1\md{243}$}, then $\nu_3(\Delta)\geq23$. As in the proof of Theorem \ref{thm1}, let $a_3=\cfrac{a}{9}$ and $u=\cfrac{-b}{8a_3}$. Since $3$ does not divide $8a_3$, then $u\in\mathbb{Z}_3$. Let $\phi=x-u$, then $F(x)=\phi^9+9u\phi^8+36u^2\phi^7+84u^3\phi^6+126u^4\phi^5+126u^5\phi^4+84u^6\phi^3+36u^7\phi^2+A\phi+B$ with $A=a+9u^8=\cfrac{\Delta}{2^{24}a^8}$ and $B=au+b+u^9=\cfrac{-b\Delta}{2^{27}a^9}$. Thus $\nu_3(A)=\nu_3(\Delta)-16$ and $\nu_3(B)=\nu_3(\Delta)-18$. It follows that $N_{\phi}(F)=S_1+S_2+S_3$ has three sides joining $(0,\nu_3(\Delta)-18)$, $(2,2)$, $(3,1)$, and $(9,0)$ with $d(S_2)=d(S_3)=1$.\\
				If $\nu_3(\Delta)$ is odd, then $d(S_1)=1$ and so $3\mathbb{Z}_K=\mathfrak{p}_1^2\mathfrak{p}_2\mathfrak{p}_3^6$ with residue degree $1$ each {ideal factor}. Hence $\nu_3(i(K))=0$.\\
				If $\nu_3(\Delta)$ is even, then $d(S_1)=2$ with ${R_{1_1}(F)(y)}=y^2+B_3\in\mathbb{F}_{\phi}[y]$. Since $(36u^7)_3\equiv 1\md{3}$, then two cases arise:
				\item[-] If $\Delta_3\not\equiv a_3\md{3}$; $\Delta_3\equiv 1\md{3}$, then {$R_{1_1}(F)(y)$} is irreducible over $\mathbb{F}_{\phi}$ and so $3\mathbb{Z}_K=\mathfrak{p}_1^2\mathfrak{p}_2\mathfrak{p}_3^6$ with $f_1=2$ and $f_2=f_3=1$. Hence $\nu_3(i(K))=0$.
				
				\item[-] If $\Delta_3\equiv a_3\md{3}$; $\Delta_3\equiv -1\md{3}$, then ${R_{1_1}(F)(y)}=(y-1)(y+1)\in\mathbb{F}_{\phi}[y]$ and so $3\mathbb{Z}_K=\mathfrak{p}_{11}\mathfrak{p}_{12}\mathfrak{p}_{21}\mathfrak{p}_{31}^6$ with residue degree $1$ each {ideal factor}. Hence $3$ divides $i(K)$. Applying \cite[Corollary, p: 230]{En}, we get $\nu_3(i(K))=1$.
			\end{enumerate}
			\item [(vii)] If $(a,b)\equiv (45,35)\md{81}$, then we have the following sub-cases:
			\begin{enumerate}
				\item[(a)] If {$a+b\equiv 80\md{243}$}, then $N_{\phi}(F)=S_1+S_2$ has two sides joining $(0,4)$, $(3,1)$, and $(9,0)$ with $d(S_2)=1$ and ${R_{1_1}(F)(y)}=y^3+y^2-y+1$ which is irreducible over $\mathbb{F}_{\phi}$. Thus $3\mathbb{Z}_K=\mathfrak{p}_1\mathfrak{p}_2^6$ with $f_{1}=3$ and $f_{2}=1$. Hence $\nu_3(i(K))=0$.
				
				\item[(b)] If {$a+b\equiv 161\md{243}$}, then $\nu_3(\Delta)\geq23$. As in the case (iv$_c$) above, let $a_3=\cfrac{a}{9}$, $u=\cfrac{-b}{8a_3}\in\mathbb{Z}_3$, and $\phi=x-u$, then $3$ divides $i(K)$ if and only if $\nu_3(\Delta)$ is even and $\Delta_3\equiv -1\md{3}$. In this case also, we have $3\mathbb{Z}_K=\mathfrak{p}_{11}\mathfrak{p}_{12}\mathfrak{p}_{21}\mathfrak{p}_{31}^6$ with residue degree $1$ each {ideal factor}. Applying \cite[Corollary, p: 230]{En}, we get $\nu_3(i(K))=1$.
				
				\item[(c)] If {$a+b\equiv -1\md{243}$}, then $N_{\phi}(F)=S_1+S_2+S_3$ has three sides joining $(0,w)$, $(1,3)$, $(3,1)$, and $(9,0)$ with $w\geq5$, $d(S_1)=d(S_3)=1$, and $R_{1_2}(F)(y)=y^2+y-1$ which is irreducible over $\mathbb{F}_{\phi}$. Thus $3\mathbb{Z}_K=\mathfrak{p}_1\mathfrak{p}_2\mathfrak{p}_3^6$ with $f_1=f_3=1$ and $f_2=2$. Hence $\nu_3(i(K))=0$.
				
			\end{enumerate}
			\item [(viii)] If $(a,b)\equiv (72,8)\md{81}$, then we have the following sub-cases:
			\begin{enumerate}
				\item[(a)] If {$a+b\equiv 80\md{243}$}, then $N_{\phi}(F)=S_1+S_2$ has two sides joining $(0,4)$, $(3,1)$, and $(9,0)$ with $d(S_2)=1$ and ${R_{1_1}(F)(y)}=y^3+y^2+1=(y-1)(y^2-y-1)\in\mathbb{F}_{\phi}[y]$. Thus $3\mathbb{Z}_K=\mathfrak{p}_{11}\mathfrak{p}_{12}\mathfrak{p}_{21}^6$ with $f_{11}=f_{21}=1$ and $f_{12}=2$. Hence $\nu_3(i(K))=0$.
				
				\item[(b)] If {$a+b\equiv 161\md{243}$}, then $N_{\phi}(F)=S_1+S_2$ has two sides joining $(0,4)$, $(3,1)$, and $(9,0)$ with $d(S_2)=1$ and ${R_{1_1}(F)(y)}=y^3+y^2-1$ which is irreducible over $\mathbb{F}_{\phi}$. Thus $3\mathbb{Z}_K=\mathfrak{p}_1\mathfrak{p}_2^6$ with $f_1=3$ and $f_2=1$. Hence $\nu_3(i(K))=0$.
				
				\item[(c)] If {$a+b\equiv -1\md{243}$}, then $\nu_3(\Delta)\geq23$. As in the case (iv$_c$) above, let $a_3=\cfrac{a}{9}$, $u=\cfrac{-b}{8a_3}\in\mathbb{Z}_3$, and $\phi=x-u$, then $3$ divides $i(K)$ if and only if $\nu_3(\Delta)$ is even and $\Delta_3\equiv -1\md{3}$. In this case also, we have $3\mathbb{Z}_K=\mathfrak{p}_{11}\mathfrak{p}_{12}\mathfrak{p}_{21}\mathfrak{p}_{31}^6$ with residue degree $1$ each {ideal factor}. Applying \cite[Corollary, p: 230]{En}, we get $\nu_3(i(K))=1$. 
			\end{enumerate}
		\end{enumerate}

		\item If $b\equiv 1\md{3}$, then $F(x)\equiv {(x+1)^9}\md{3}$. Let $\phi=x+1$. {Then} $F(x)=\phi^9-9\phi^8+36\phi^7-84\phi^6+126\phi^5-126\phi^4+84\phi^3-36\phi^2+(a+9)\phi-a+b-1$.
		\begin{enumerate}
			\item[(i)] If $(a,b)\in\{(0,4),(0,7),(3,1),(3,7),(6,1),(6,4)\}\md{9}$, then by Theorem \ref{thm1}, $\nu_3(i(K))=0$.
			
			\item[(ii)] {The treatment of the other cases is similar to the case $b\equiv -1\md{3}$ above. Table $9$ summarizes the obtained results.}
				\begin{table}[htbp]
				\centering
				\caption{}
				\begin{tabular}{|l|l|c|c|c|}
					\hline
					\multicolumn{2}{|c|}{Cases}&$3\mathbb{Z}_K$&\multicolumn{1}{|c|}{$f_i$}&$\nu_3(i(K))$\\
					\hline
					\multicolumn{2}{|l|}{$(a,b)\in\{(3,4),(6,7)\}\md{9}$}&$\mathfrak{p}_1\mathfrak{p}_2^8$&$f_i=1$&\multirow{10}{*}{$0$}\\
					\cline{1-4}
					\multicolumn{2}{|l|}{$(a,b)\in\{(0,10),(0,19),(9,1),(9,19),(18,1),(18,10)\}$}&\multirow{2}{*}{$\mathfrak{p}_1^3\mathfrak{p}_2^6$}&\multirow{2}{*}{$f_i=1$}&\\
					\multicolumn{2}{|c|}{$\md{27}$}&&&\\
					\cline{1-4}
					\multicolumn{2}{|l|}{$(a,b)\in\{(0,1),(9,10)\}\md{27}$}&$\mathfrak{p}_1\mathfrak{p}_2^2\mathfrak{p}_3^6$& $f_i=1$&\\
					\cline{1-4} 
					\multicolumn{2}{|l|}{$(a,b)\in\{(18,46),(18,73),(45,19),(45,73),(72,19),$}&\multirow{2}{*}{$\mathfrak{p}_1^3\mathfrak{p}_2^6$}&\multirow{2}{*}{$f_i=1$}&\\
					\multicolumn{2}{|c|}{$(72,46)\}\md{81}$}&&&\\
					\cline{1-4}
					\multicolumn{2}{|l|}{$(a,b)\equiv (18,19)\md{81}$ and $b-a\equiv 82\md{243}$}&$\mathfrak{p}_{1}\mathfrak{p}_2^6$&$f_1=3$, $f_2=1$&\\
					\cline{1-4} 
					\multicolumn{2}{|l|}{$(a,b)\equiv (18,19)\md{81}$ and $b-a\equiv 163\md{243}$}&$\mathfrak{p}_{1}\mathfrak{p}_{2}\mathfrak{p}_{3}^6$&$f_{1}=f_{3}=1$, $f_{2}=2$&\\
					\cline{1-4}
					$(a,b)\equiv (18,19)\md{81}$ &$\nu_3(\Delta)$ is odd&$\mathfrak{p}_1^2\mathfrak{p}_2\mathfrak{p}_3^6$&$f_i=1$&\\
					\cline{2-4}
					\multicolumn{1}{|c|}{and} &$\nu_3(\Delta)$ is even, $\Delta_3\equiv 1\md{3}$&$\mathfrak{p}_1^2\mathfrak{p}_2\mathfrak{p}_3^6$&$f_1=2$, $f_2=f_3=1$&\\
					\cline{2-5}
					$b-a\equiv 1\md{243}$&$\nu_3(\Delta)$ is even, $\Delta_3\equiv -1\md{3}$&$\mathfrak{p}_{1}\mathfrak{p}_{2}\mathfrak{p}_{3}\mathfrak{p}_{4}^6$&$f_i=1$&$1$\\
					\hline
					\multicolumn{2}{|l|}{$(a,b)\equiv (45,46)\md{81}$ and $b-a\equiv 163\md{243}$}&$\mathfrak{p}_1\mathfrak{p}_2^6$& $f_{1}=3$, $f_{2}=1$&\multirow{4}{*}{$0$}\\
					\cline{1-4}
					\multicolumn{2}{|l|}{$(a,b)\equiv (45,46)\md{81}$ and $b-a\equiv 1\md{243}$}&$\mathfrak{p}_1\mathfrak{p}_2\mathfrak{p}_3^6$&$f_1=f_3=1$, $f_2=2$&\\
					\cline{1-4}
					$(a,b)\equiv (45,46)\md{81}$ &$\nu_3(\Delta)$ is odd&$\mathfrak{p}_1^2\mathfrak{p}_2\mathfrak{p}_3^6$&$f_i=1$&\\
					\cline{2-4}
					\multicolumn{1}{|c|}{and} &$\nu_3(\Delta)$ is even, $\Delta_3\equiv 1\md{3}$&$\mathfrak{p}_1^2\mathfrak{p}_2\mathfrak{p}_3^6$&$f_1=2$, $f_2=f_3=1$&\\
					\cline{2-5}
					$b-a\equiv 82\md{243}$&$\nu_3(\Delta)$ is even, $\Delta_3\equiv -1\md{3}$&$\mathfrak{p}_{1}\mathfrak{p}_{2}\mathfrak{p}_{3}\mathfrak{p}_{4}^6$&$f_i=1$&$1$\\
					\hline
					\multicolumn{2}{|l|}{$(a,b)\equiv (72,73)\md{81}$ and $b-a\equiv 82\md{243}$}&$\mathfrak{p}_1\mathfrak{p}_2^6$& $f_{1}=3$, $f_{2}=1$&\multirow{4}{*}{$0$}\\
					\cline{1-4}
					\multicolumn{2}{|l|}{$(a,b)\equiv (72,73)\md{81}$ and $b-a\equiv 163\md{243}$}&$\mathfrak{p}_1\mathfrak{p}_2\mathfrak{p}_3^6$&$f_1=f_3=1$, $f_2=2$&\\
					\cline{1-4}
					
					$(a,b)\equiv (72,73)\md{81}$ &$\nu_3(\Delta)$ is odd&$\mathfrak{p}_1^2\mathfrak{p}_2\mathfrak{p}_3^6$&$f_i=1$&\\
					\cline{2-4}
					\multicolumn{1}{|c|}{and} &$\nu_3(\Delta)$ is even, $\Delta_3\equiv 1\md{3}$&$\mathfrak{p}_1^2\mathfrak{p}_2\mathfrak{p}_3^6$&$f_1=2$, $f_2=f_3=1$&\\
					\cline{2-5}
					$b-a\equiv 1\md{243}$&$\nu_3(\Delta)$ is even, $\Delta_3\equiv -1\md{3}$&$\mathfrak{p}_{1}\mathfrak{p}_{2}\mathfrak{p}_{3}\mathfrak{p}_{4}^6$&$f_i=1$&$1$\\
					\hline
				\end{tabular}
			\end{table}
		\end{enumerate}
		\item If $b\equiv 0\md{3}$, then $F(x)\equiv x^9\md{3}$. Let $\phi=x$.  {Then} $F(x)=\phi^9+a\phi+b$. By assumption, $\nu_3(b)\leq8$ or $\nu_3(a)\leq7$.\\
		If $8\nu_3(b)<9\nu_3(a)$, then $N_{\phi}(F)=S_1$ has a single side joining $(0,\nu_3(b))$ and $(9,0)$ with degree $d\in\{1,3\}$. 
		\begin{enumerate}
			\item[(i)] If $d=1$, then {$R_{1_1}(F)(y)$} is irreducible as it is of degree $1$. Thus $3\mathbb{Z}_K=\mathfrak{p}_1^9$ with residue degree $1$. Hence $\nu_3(i(K))=0$.
			
			\item[(ii)] If $d=3$, then $-\lambda_1=-1/3,-2/3$ is the slope of $S_1$. Since $e_1=3$ divides the the ramification index of any prime ideal of $\mathbb{Z}_K$ lying above $3$, then there is at most three prime ideals of $\mathbb{Z}_K$ lying above $3$ with residue degree $1$ each {ideal factor}. Hence $\nu_3(i(K))=0$.
		\end{enumerate}
		
		If $8\nu_3(b)>9\nu_3(a)$, then {$N_{\phi}(F)=S_1+S_2$ has two sides joining $(0,\nu_3(b))$, $(1,\nu_3(a))$, and $(9,0)$ with $d(S_1)=1$ and $d(S_2)\in\{1,2,4\}$ since $\nu_3(a)\leq 7$. Let $d=d(S_2)$. Then we have the following cases:} 
		\begin{enumerate}
			\item[(i)] {If $d=1$; $\nu_3(a)\in\{1,3,5,7\}$, then $3\mathbb{Z}_K=\mathfrak{p}_1\mathfrak{p}_2^8$ with residue degree $1$ each ideal factor. Hence $\nu_3(i(K))=0$}.
			
			\item[(ii)] {If $d=2$, then ${R_{1_2}(F)(y)}=a_3y^2+b_3\in\mathbb{F}_{\phi}[y]$; that is ${R_{1_2}(F)(y)}=\pm(y^2+1)$ or  ${R_{1_2}(F)(y)}=\pm(y-1)(y+1)$. Thus $3\mathbb{Z}_K=\mathfrak{p}_1\mathfrak{p}_2^4$ with $f_1=1$ and $f_2=2$  or  $3\mathbb{Z}_K=\mathfrak{p}_{11}\mathfrak{p}_{21}^4\mathfrak{p}_{22}^4$ with $f_{11}=f_{21}=f_{22}=1$ respectively. Hence $\nu_3(i(K))=0$}.
			
			\item[(iii)] {If $d=4$, then ${R_{1_2}(F)(y)}=a_3y^4+b_3\in\mathbb{F}_{\phi}[y]$; that is $R_{1_2}(F)(y)=\pm(y^4+1)=\pm(y^2-y-1)(y^2+y-1)$ or  ${R_{1_2}(F)(y)}=\pm(y^4-1)=\pm(y-1)(y+1)(y^2+1)$. Thus $3\mathbb{Z}_K=\mathfrak{p}_{11}\mathfrak{p}_{21}^2\mathfrak{p}_{22}^2$ with $f_{11}=1$ and $f_{21}=f_{22}=2$ or  $3\mathbb{Z}_K=\mathfrak{p}_{11}\mathfrak{p}_{21}^2\mathfrak{p}_{22}^2\mathfrak{p}_{23}^2$ with $f_{11}=f_{21}=f_{22}=1$ and $f_{23}=2$ respectively. Hence $\nu_3(i(K))=0$}.
			
		\end{enumerate}

	\end{enumerate}
	\begin{flushright}
		$\square$
	\end{flushright}
	
	\textit{\textbf{Proof of Theorem \ref{thmp}}}.\\
	For $p=5$. Since $\Delta=2^{24}a^{9}+3^{18}b^{8}$, {$\nu_5(\Delta)\geq1$ if and only if $(a,b)\in\{(0,0),(1,1),(1,2),$\\$(1,3),(1,4)\}\md{5}$. Thanks to the index formula $(1.1)$, $5$ can divides the index $i(K)$ only if $(a,b)\in\{(0,0),(1,1),(1,2),(1,3),(1,4)\}\md{5}$.}
	\begin{enumerate}
		\item For $(a,b)\in\{(1,1),(1,2),(1,3),(1,4)\}\md{5}$, one can easily check that $F(x)\equiv \phi_{i1}\cdot\phi_{i2}^2\md{5}$ with $\deg(\phi_{i1})=7$, $\deg(\phi_{i2})=1$, and $\phi_{ij}$ is irreducible over $\mathbb{F}_5$ for every $i=1,\dots,4$ and $j=1,2$. Thus there is at most two prime ideals of $\mathbb{Z}_K$ lying above $5$ with residue degree $1$ each {ideal factor}. Hence $\nu_5(i(K))=0$.
		
		\item If $(a,b)\equiv (0,0)\md{5}$, then $F(x)\equiv x^{9}\md{5}$. Let $\phi=x$. Then $F(x)=\phi^{9}+a\phi+b$.
		\begin{enumerate}
			\item[(i)] If $8\nu_5(b)<9\nu_5(a)$, then $N_{\phi}(F)=S_1$ has a single side joining $(0,\nu_5(b))$ and $(9,0)$ with degree $d\in\{1,3\}$.
			\begin{enumerate}
				\item[(a)] If $d=1$, then {$R_{1_1}(F)(y)$} is irreducible as it is of degree $1$. Thus $5\mathbb{Z}_K=\mathfrak{p}_1^9$ with residue degree $1$. Hence $\nu_5(i(K))=0$.
				
				\item[(b)] If $d=3$, then ${R_{1_1}(F)(y)}=y^3+b_5\in\mathbb{F}_{\phi}[y]$. One can check easily that, for every value $b_5\in\mathbb{F}_{\phi}^*$, $R_{1_1}(F)(y)=\psi_1\cdot\psi_2$ with $\deg(\psi_1)=1$, $\deg(\psi_2)=2$, and $\psi_i$ is irreducible over $\mathbb{F}_{\phi}$. Thus $5\mathbb{Z}_K=\mathfrak{p}_1^3\mathfrak{p}_2^3$ with $f_1=1$ and $f_2=2$. Hence $\nu_5(i(K))=0$.
			\end{enumerate}
			\item[(ii)] If $8\nu_5(b)>9\nu_5(a)$, then {$N_{\phi}(F)=S_1+S_2$ has two sides joining $(0,\nu_5(b))$, $(1,\nu_5(a))$, and $(9,0)$ with $d(S_1)=1$ and $d(S_2)\in\{1,2,4\}$ since $\nu_5(a)\leq 7$. Thus $\phi$ can provides at most five prime ideal of $\mathbb{Z}_K$ lying above $5$ with residue degree $1$ each ideal factor. Hence $\nu_5(i(K))=0$.}
		\end{enumerate}
	\end{enumerate}
	For $p=7$. Since $\Delta=2^{24}a^{9}+3^{18}b^{8}$, $\nu_7(\Delta)\geq1$ if and only if $(a,b)\in\{(0,0),(3,1),(3,6),$\\$(5,1),(5,6),(6,1),(6,6)\}\md{7}$. Thanks to the index formula $(1.1)$, $7$ can divides the index $i(K)$ only if $(a,b)\in\{(0,0),(3,1),(3,6),(5,1),(5,6),(6,1),(6,6)\}\md{7}$.
	\begin{enumerate}
		\item For $(a,b)\in\{(3,1),(3,6),(5,1),(5,6),(6,1),(6,6)\}\md{7}$, one can easily check that $F(x)\equiv \phi_{i1}\cdot\phi_{i2}\cdot\phi_{i3}^2\md{7}$ with $\deg(\phi_{i1})=4$, $\deg(\phi_{i2})=3$, $\deg(\phi_{i2})=1$, and $\phi_{ij}$ is irreducible over $\mathbb{F}_7$ for every $i=1,\dots,6$ and $j=1,2,3$. Thus there is at most two prime ideals of $\mathbb{Z}_K$ lying above $7$ with residue degree $1$ each {ideal factor}. Hence $\nu_7(i(K))=0$.
		
		\item If $(a,b)\equiv (0,0)\md{7}$, then $F(x)\equiv x^{9}\md{7}$. Let $\phi=x$. Then $F(x)=\phi^{9}+a\phi+b$.
		\begin{enumerate}
			\item[(i)] If $8\nu_7(b)<9\nu_7(a)$, then $N_{\phi}(F)=S_1$ has a single side joining $(0,\nu_7(b))$ and $(9,0)$ with degree $d\in\{1,3\}$. 
			\begin{enumerate}
				\item[(a)] If $d=1$, then $R_{1_1}(F)(y)$ is irreducible as it is of degree $1$. Thus $7\mathbb{Z}_K=\mathfrak{p}_1^9$ with residue degree $1$. Hence $\nu_7(i(K))=0$.
				
				\item[(b)] If $d=3$, then $R_{1_1}(F)(y)=y^3+b_5\in\mathbb{F}_{\phi}[y]$. Thus there is at most three prime ideals of $\mathbb{Z}_K$ lying above $7$ with residue degree $1$ each {ideal factor}. Hence $\nu_7(i(K))=0$.
			\end{enumerate}
			\item[(ii)] If $8\nu_7(b)>9\nu_7(a)$, then {$N_{\phi}(F)=S_1+S_2$ has two sides joining $(0,\nu_7(b))$, $(1,\nu_7(a))$, and $(9,0)$ with $d(S_1)=1$ and $d(S_2)\in\{1,2,4\}$ since $\nu_7(a)\leq 7$. Thus $\phi$ can provides at most five prime ideal of $\mathbb{Z}_K$ lying above $7$ with residue degree $1$ each ideal factor. Hence $\nu_7(i(K))=0$.}
		\end{enumerate} 
	\end{enumerate}
	For $p\geq11$, since there is at most $9$ prime ideals of $\mathbb{Z}_{K}$ lying above $p$ with residue degree $1$ each, and there is at least $p\geq 11$ monic irreducible polynomial of degree $f$ in $\mathbb{F}_p[x]$ for every positive integer $f$, we conclude that $p$ does not divide $i(K)$.
	\begin{flushright}
		$\square$
	\end{flushright}
	\section{Examples}
	Let $F(x)=x^{9}+ax+b\in \mathbb{Z}[x]$ be a monic irreducible polynomial and $K$ the nonic number field generated by a complex root of $F(x)$.
	\begin{enumerate}
		\item For $a=51$ and $b=122$, we have $(a,b)\equiv (3,2)\md{4}$, $(a,b)\equiv(6,5) \md{9}$, and for every rational prime $p\notin\{2,3\}$, $\nu_p(\Delta)\leq1$. By Theorem $\ref{thm1}$, $\mathbb{Z}[\alpha]$ is integrally closed and so $K$ is monogenic. Hence $i(K)=1$.
		
		\item For $a=35$ and $b=20$, we have $(a,b)\equiv(3,4)\md{8}$, then by Theorem $\ref{thmp2}$, $i(K)$ is even. Hence $K$ is {not monogenic.}
		\item For $a=1392$ and $b=768$, we have $(a,b)\equiv(368,256)\md{512}$, then by Theorem {$\ref{thmp2}$}, $\nu_2(i(K))=1$. On the other hand, $F(x)$ is $3$-Eisenstein, then $\nu_3(i(K))=0$. We conclude that $i(K)=2$. Hence $K$ is not monogenic.
		
		\item For $a=126$ and $b=40130$, we have {$(a,b)\equiv (45,35)\md{81}$, $a+b\equiv 161\md{243}$}, $\nu_3(\Delta)=26$, and $\Delta_3\equiv -1\md{3}$, then by Theorem $\ref{thmp3}$, $\nu_3(i(K))=1$. On the other hand, $F(x)$ is $2$-Eisenstein, then $\nu_2(i(K))=0$. We conclude that $i(K)=3$. Hence $K$ is not monogenic.
		
		\item For $a=15381$ and $b=6634$, we have $(a,b)\equiv (1,2)\md{4}$, then by Theorem $\ref{thmp2}$, $\nu_2(i(K))=1$. On the other hand, {$(a,b)\equiv (72,73)\md{81}$, $b-a\equiv 1\md{243}$}, $\nu_3(\Delta)=24$, and $\Delta_3\equiv -1\md{3}$, then by Theorem $\ref{thmp3}$, $\nu_3(i(K))=1$. We conclude that $i(K)=6$. Hence $K$ is not monogenic.
		
		\item For $a=183$ and $b=296$, we have $(a,b)\equiv (7,8)\md{16}$ and $\nu_2(\Delta)=29$, then by Theorem $\ref{thmp2}$, $\nu_2(i(K))=3$. On the other hand, $(a,b)\equiv (21,53)\md{81}$, then by Theorem $\ref{thmp3}$, $\nu_3(i(K))=0$. We conclude that $i(K)=8$. Hence $K$ is not monogenic.
		
		\item For $a=7335$ and $b=24184$, we have $(a,b)\equiv (7,8)\md{16}$, $\nu_2(\Delta)=28$, and $\Delta_2\equiv 3\md{8}$, then by Theorem $\ref{thmp2}$, $\nu_2(i(K))=3$. On the other hand, {$(a,b)\equiv (45,46)\md{243}$, $b-a\equiv82\md{243}$}, $\nu_3(\Delta)=24$, and $\Delta_3\equiv -1\md{3}$, then by Theorem $\ref{thmp3}$, $\nu_3(i(K))=1$. We conclude that $i(K)=24$. Hence $K$ is not monogenic.
	\end{enumerate}
	\text{ }\newline
	{\bf Conflict of interest}\\
	Not Applicable.\\
	{\bf Data availability}\\
	Not applicable.\\
	{\bf Author Contribution and Funding Statement}\\
	Not applicable.

\end{document}